\documentclass{amsart}

\usepackage{amssymb}
\usepackage[all]{xy}
\usepackage{hyperref}

\usepackage{enumitem}   

\setlist[enumerate]{itemsep=.2em,topsep=.2em,leftmargin=1.25em,itemindent=2.0em}

\newtheorem{thm}{Theorem}

\newtheorem{lem}[thm]{Lemma}
\newtheorem{cor}[thm]{Corollary}

\newtheorem{prop}[thm]{Proposition}

\newtheorem{conj}[thm]{Conjecture}

\theoremstyle{definition}

\newtheorem{say}[thm]{}
\newtheorem{exmp}[thm]{Example}

\newtheorem{rem}[thm]{Remark}

\newtheorem*{ack}{Acknowledgments}      
\newtheorem{notation}[thm]{Notation}   
  
\newtheorem{defn-thm}[thm]{Definition--Theorem}  
\newtheorem{defn-lem}[thm]{Definition--Lemma}  

\theoremstyle{remark}

\renewcommand{\c}[0]{{\mathbb C}}

\newcommand{\z}[0]{{\mathbb Z}}
\newcommand{\n}[0]{{\mathbb N}}

\renewcommand{\a}[0]{{\mathbb A}}

\newcommand{\p}[0]{{\mathbb P}}

\newcommand{\qtq}[1]{\quad\mbox{#1}\quad}

\newcommand{\aut}[0]{\operatorname{Aut}}

\newcommand{\rdown}[1]{\lfloor{#1}\rfloor}

\newcommand{\tsum}[0]{\textstyle{\sum}}

\def\loccoh#1.#2.#3.#4.{H^{#1}_{#2}(#3,#4)}

\DeclareMathAlphabet{\mathchanc}{OT1}{pzc}%
                                {m}{it}

\usepackage[all]{xy}\xyoption{dvips}

\begin{document}
\bibliographystyle{amsalpha}


 \title[Mordell-Schinzel conjecture]{The Mordell-Schinzel conjecture \\ for cubic diophantine equations}
 \author{J\'anos Koll\'ar and Jennifer Li}

 \begin{abstract}
   Building on the works of Mordell (1952) and Schinzel (2015), we prove that   cubic   diophantine equations
 $xyz=G(x,y)$ have infinitely many solutions.
       \end{abstract}

               \maketitle

 Mordell stated \cite{MR51852} that  the diophantine equation
 $xyz=G(x,y)$ has infinitely many solutions for every polynomial $G\in \z[x,y]$.
 This turned out to be not completely true; counterexamples for $\deg G=2$ were already known to  Jacobs\-thal    \cite{MR1557037}.
 A detailed study was undertaken by Schinzel \cite{MR3320468, MR3826636}; see also \cite[Sec.5.5]{MR4823864}.

 \begin{notation}\label{not.1}
 By changing $z$, the terms in $G$ that are divisible by $xy$ can be absorbed into $xyz$, thus, setting $c:=G(0,0)$ we can simplify the equation to 
 $$
 xyz= A(x)+B(y) -c,
 \eqno{(\ref{not.1}.1)}
 $$
 where $A(x):=\tsum_{i=0}^n a_i x^i$,
 $B(y):= \tsum_{j=0}^m b_j y^j$, and $c=a_0=b_0$.
 (Thus $c$ appears  3 times in (\ref{not.1}.1).)
  We write $a:=a_n$, $b:=b_m$ and, as in  \cite{MR3320468}, 
  assume from now on that  $abc\neq 0$.
 The corresponding surfaces are denoted by
 $$
 S_{A,B}:=\bigl(xyz=A(x)+B(y) -c\bigr)\subset \a^3.
  \eqno{(\ref{not.1}.2)}
 $$
 As in \cite{k-vp-1}, these surfaces can be defined over an arbitrary integral domain $R$. Here we work over $\z$, but some general results are stated over $\c$.

 Writing the equation as 
 $$
 xyz= ax^n+by^m+c + \tsum_{i=1}^{n-1} a_i x^i + \tsum_{j=1}^{m-1} b_j y^j
 \eqno{(\ref{not.1}.3)}
 $$
 suggests that  one should view them as perturbations of the
 trinomial series $xyz= ax^n+by^m+c$.
 See Paragraph~\ref{geom.say} for a geometric viewpoint.

 \end{notation}

 Somewhat surprisingly, the cases where $\min\{m, n\}\in \{1,2\}$ are especially problematic. In particular, the case  $m=n=2$ is not completely settled
 \cite{MR3320468, MR3826636}. Several of these equations are also discussed in \cite[Chap.5]{MR4823864}.

 However, the situation seems better for $m,n\geq 3$. Based on
 \cite[Thm.4]{MR3320468}, we propose the following conjecture.

 \begin{conj}[Mordell-Schinzel]\label{m-s-conj}  If $m, n\geq 3$ and $A(t), B(t)\in \z[t]$,  then
   (\ref{not.1}.1) has infinitely many integer solutions.
 \end{conj}

 For historical accuracy, note that neither author stated this as a conjecture.  \cite{MR51852} asserted this as a true statement, while
 \cite[Thm.4]{MR3320468} proved most cases, without speculating about the general version; see Paragraph~\ref{ms.comments} for further details.  To the best of our knowledge, there is not a single
 pair $(m,n)$ for which Conjecture~\ref{m-s-conj} has been fully proven.
 Our aim is to settle the simplest case.

 \begin{thm}\label{m-s-3.thm}  For every $a,b\in\z\setminus\{0\}$ and $a_1, a_2, b_1, b_2, c\in \z$, the  cubic surface
   $$
   \bigl(xyz=ax^3+by^3+c+a_2x^2+a_1x+b_2y^2+b_1y\bigr)\subset
   \a^3
   $$
   has  infinitely many integral points.
 \end{thm}

 \begin{say}[Plan of the proof] \label{m-s-3.thm.plan}
   If $c=0$ then we have the trivial solutions $(0,0,z)$, so we may  assume that $c\neq 0$.
   
     \cite[Thm.4]{MR3320468} establishes all the cases when
   $|abc|>1$.  Thus it remains to deal with the possibilities 
   $a,b, c\in \{1, -1\}$.

These are very special, but for $m=n=2$  all
known equations  (\ref{not.1}.3) with only finitely many solutions satisfy $a,b,c\in \{1, -1\}$;
see \cite{MR1557037, MR3826636}.  
So there is reason to think that the  $a,b, c\in \{1, -1\}$ cases could behave differently.

   By changing the signs of $x,y,z$, we can achieve $a=b=c=1$. Thus  from now on we work with the surfaces
$$
 S_{a_1, a_2, b_1, b_2}:=  \bigl(xyz=x^3+y^3+1+a_2x^2+a_1x+b_2y^2+b_1y\bigr)\subset
 \a^3.
 \eqno{(\ref{m-s-3.thm.plan}.1)}
   $$
 The observation of \cite{k-vp-1} is that these surfaces have an infinite automorphism group over $\z$.
 In retrospect, these automorphisms  are behind the computations in
\cite{MR51852}; more explicit versions appear in  \cite{MR0769779, MR3320468, bog-sha}. 

 Thus we aim to find {\em one} solution,  and hope to prove that its orbit under the automorphism group is infinite.
The easiest is to  start with the {\it trivial solutions} $(x_0, y_0, z_0)$, where $x_0=\pm 1, y_0=\pm 1$.

In most cases the automorphism group is infinite cyclic, but its
 generator is unexpectedly complicated. It is  given by
 polynomials of degrees $13$,  $34$, and  $55$,   containing $110$, $998$, and  $2881$   monomials; see \cite[Exmp.2]{k-vp-1}.  So it is almost impossible to compute with.

 Therefore we  work  instead with a closely related  isomorphism-groupoid
 $\Sigma_{A,B}$ involving 4 surfaces, called the {\it companion surfaces/equations} of (\ref{not.1}.2); see Paragraph~\ref{g.isom.say}.

 Using the explicit forms of the generators of $\Sigma_{A,B}$
  it is easy to see that if $\max\{|x_0|, |y_0|\}$ is large enough, then
 the   $\Sigma_{A,B}$-orbit of $p_0=(x_0, y_0, z_0)$ is infinite; see
 Corollary~\ref{inf.orb.cor.1}. This holds for all $m, n\geq 3$.

 The rest of the argument assumes that $m=n=3$.
 The required condition on  $\max\{|x_0|, |y_0|\}$ is not satisfied by the trivial solutions of (\ref{m-s-3.thm.plan}.1), but straightforward estimates show that if
 $\max\{|a_1|+|a_2|, |b_1|+|b_2|\}\geq 7$,
 then there is a trival solution  of one of the companion equations,
whose $\Sigma_{A,B}$-orbit is infinite.
This is proved in 
 Corollary~\ref{7.sols.cor}.

This shows that  Conjecture~\ref{m-s-conj} holds for cubics, 
save possibly for a finite list of exceptions, all of which satisfy
$|a_1|+|a_2|, |b_1|+ |b_2|\leq 6$. 

For these equations we rely on computer searches.
First, in Proposition~\ref{fin.triv.orb.prop} we list all cubics
$ S_{a_1, a_2, b_1, b_2}$ for which the 
  $\Sigma_{A,B}$-orbits of all trivial solutions are finite.
There are only 3 such, up to isomorphism.

In these cases we searched for other solutions. Once we find one
integral point  $(x_0, y_0, z_0)$ with $\max\{|x_0|, |y_0|\}\geq 6$,
its $\Sigma_{A,B}$-orbit is  infinite  by 
Corollary~\ref{inf.orb.cor.1}. \qed
\end{say}

 \begin{exmp}
The above  approach does not work directly for
$m=n=4$, since there are infinitely many quartic polynomials
$A(x)+B(y)-c$   with $|abc|=1$ for which
the 
$\Sigma_{A,B}$-orbits of all trivial solutions are finite.

   If $A(x)$ and $B(y)$  map  $\{1, -1\}$ to itself, then
  $\Sigma_{A,B}$  maps trivial solutions to  trivial solutions.
  Thus the orbit of all trivial solutions is finite.
  An infinite series of such polynomials is given by
  $x^4-x^2+1+r(x^3-x)$ for $r\in \z$.
   \end{exmp}

 \begin{rem} \label{m.m-s-3.thm.plan} As  in  \cite{bog-sha}, we also get some information about the equations 
   $$
 mxyz=x^3+y^3+1+a_2x^2+a_1x+b_2y^2+b_1y
 \eqno{(\ref{m.m-s-3.thm.plan}.1)}
   $$
 for fixed $m\in \z\setminus\{0\}$.
 
 For $m\neq \pm 1$, these equations do not have automorphisms,
 so we look at the action  of the automorphism group on the solutions of
 (\ref{m-s-3.thm.plan}.1) modulo $m$.
 If (\ref{m-s-3.thm.plan}.1) has a solution
 $(x_0, y_0, z_0)$ whose orbit is infinite, and
 $z_0\equiv 0\mod m$, then the orbit contains infinitely many other solutions 
 $(x_r, y_r, z_r)$ such that  $z_r\equiv 0\mod m$.
 All these give solutions of (\ref{m.m-s-3.thm.plan}.1). 
 
Thus, using Corollary~\ref{inf.orb.cor.1},
we get that if there is at least one solution  $(x_0, y_0, z_0)$ of (\ref{m.m-s-3.thm.plan}.1) satisfying
 $$
 \max\{|x_0|, |y_0|\}> 2+
 \max\{|a_2|+|a_1|, |b_2|+|b_1|\},
 $$
 then there are infinitely many.
 \end{rem}

 \begin{say}[The isomorphisms $\sigma_x, \sigma_y$]\label{s.isom.say}
   Let $S_{A, B}$ be as in (\ref{not.1}.2--3) and 
   assume that  $a=b=c=1$. Set
   $$
   \bar A(x):=x^nA(x^{-1}) \qtq{and} \bar B(y):=y^mB(y^{-1}).
   \eqno{(\ref{s.isom.say}.1)}
   $$
   The surfaces  $ S_{A,B},  S_{\bar A,  B}, S_{ A, \bar B}, S_{\bar A, \bar B}$ are the {\it companion surfaces} of
 $ S_{A,B}$; their equations are the 
   {\it   companion equations.}

   By  \cite[27]{k-vp-1} there are isomorphisms
   $$
   \sigma_{y,A,B}: S_{A,B}\cong S_{\bar A,B}
   \qtq{and} \sigma_{x,A,B}: S_{A,B}\cong S_{ A,\bar B}.
   \eqno{(\ref{s.isom.say}.2)}
     $$
    The most transparent are the rational forms on the $(x,y)$-coordinates
   $$
   \sigma_{y,A,B}: (x, y)\mapsto \bigl(B(y)/x, y\bigr)
   \qtq{and} \sigma_{x,A,B}:(x, y)\mapsto \bigl(x, A(x)/y\bigr);
   \eqno{(\ref{s.isom.say}.3)}
     $$
     see also  \cite[Formulas~52-53]{MR3320468}.
     For  $xy=0$  we need to use the  versions
$$
     (x, y)\mapsto \bigl(yz-\tsum_{i=1}^na_i x^{i-1}, y\bigr) \qtq{and}
      (x, y)\mapsto \bigl(x, xz-\tsum_{j=1}^mb_j y^{j-1}\bigr).
   \eqno{(\ref{s.isom.say}.4)}
     $$

   These are involutions in the sense that $\sigma_{y,\bar A,B}\circ \sigma_{y,A,B}=1$
   and $\sigma_{x, A,\bar B}\circ \sigma_{x,A,B}=1$.

     We will usually write  $\sigma_{x}:=\sigma_{x,A,B}$ and
   $\sigma_{y}:=\sigma_{y,A,B}$. Then   $\sigma_{x}^2=1$ and $\sigma_{y}^2=1$.

     \medskip

     {\it Remark \ref{s.isom.say}.5.}  If $(0, y_0, z_0)$ is a point on $S_{A,B}$ then  $y_0\neq 0$ since $c\neq 0$, and $z_0$ can be arbitrary.
     Thus applying $\sigma_y$ as in  (\ref{s.isom.say}.4) we get points with
     $x_1$ arbitrarily large. Thus Conjecture~\ref{m-s-conj} holds for such surfaces. We can thus always work with the simpler forms (\ref{s.isom.say}.3) of  $\sigma_x$ and  $\sigma_y$.
 \end{say}

 \begin{say}[The isomorphism groupoid]\label{g.isom.say}
   Let $S_{A, B}$ be as in (\ref{not.1}.2--3) and 
   assume that  $a=b=c=1$.
   The isomorphisms  (\ref{s.isom.say}.2) of the  companion surfaces
   sit in a diagram
   $$
\begin{array}{ccc}
  S_{A,B} & \stackrel{\sigma_x}{\longleftrightarrow} &  S_{A, \bar B}\\[1ex]
  \sigma_y\updownarrow \hphantom{\sigma_y} &&  \hphantom{\sigma_y}\updownarrow\sigma_y \\
  S_{\bar A, B} & \stackrel{\sigma_x}{\longleftrightarrow} &  S_{\bar A, \bar B}
\end{array}
 \eqno{(\ref{g.isom.say}.1)}
 $$
 The  groupoid generated by these 8 isomorphisms $\sigma_x, \sigma_y$
 is denoted by  $\Sigma_{A,B}$.

By \cite[Thm.3]{k-vp-1},  going around the diagram (\ref{g.isom.say}.1)
   clockwise gives an infinite order automorphism
   $$
   \sigma_{A,B}:=\sigma_y\circ\sigma_x\circ\sigma_y\circ\sigma_x\in
   \aut_{\z}(S_{A,B}). 
   \eqno{(\ref{g.isom.say}.2)}
   $$
Moreover, if $m=n=3$ then $\langle \sigma_{A,B}\rangle$ has finite index in
 $\aut_{\c}(S_{A,B})$  by \cite[Thm.29]{k-vp-1}, and
 every infinite $\sigma_{A,B}$-orbit is 
 Zariski dense  by \cite[51]{k-vp-1}. Most likely these also hold for all $n, m\geq 3$. 
 
 Our plan is to start with an integral point on one of the  companion surfaces of
 $ S_{A,B}$,  and then move it around using
  $\Sigma_{A,B}$
 to get infinitely many other integral points.
 We need an  elementary root estimate.

 \end{say}

\begin{lem} \label{elem.pos.lem} Let $f(x)=x^n+\tsum_{i=0}^{n-1} a_i x^i$ be a complex polynomial.
  Then
  \begin{enumerate}
  \item all roots $r_j$ of $f(x)$ satisfy $|r_j|\leq \max\bigl\{1, \tsum_{i=0}^{n-1} |a_i|\bigr\}$, and
  \item  $|f(x)|>|x^m|$ for $0\leq m<n$ and
    $|x|\leq -1+\tsum_{i=0}^{n-1} |a_i|$.   \qed
    \end{enumerate}
\end{lem}

We measure the size of points using the norm
 $|(x,y,z)|:=\max\{|x|, |y|\}$.

\begin{cor}\label{inf.orb.cor.1}
With $S_{A, B}$  as in (\ref{not.1}.2--3),   
   assume that  $a=b=c=1$ and $n, m\geq 3$. 
  Let $p_0=(x_0, y_0, z_0)\in S_{A,B}$ be a complex point. 
  If
  $$
  |p_0|>  1+\max\bigl\{1, \tsum_{i=0}^{n-1} |a_i|,   \tsum_{j=0}^{m-1} |b_j|\bigr\},
  $$
  then
  $
  \max\bigl\{ |\sigma_x(p_0)|, |\sigma_y(p_0)|\bigr\} > |p_0|.
  $
Therefore 
   the  $\Sigma_{A,B}$-orbit of $p_0$ is infinite.
\end{cor}

Proof.  We may assume that $|x_0|\geq |y_0|$ and consider
$$
p_1:=\sigma_x(p_0)=\bigl(x_0, A(x_0)/y_0, z_1\bigr).
$$
Here $\bigl|A(x_0)\bigr|>|x_0|^2$ by (\ref{elem.pos.lem}.2) and  $y_0\neq 0$ since
$x_0$ is not a root of $A(x)$. So
$\bigl|A(x_0)/y_0\bigr|>|x_0|$.

Alternating $\sigma_x, \sigma_y$ gives a sequence of points
$p_i$ in the   $\Sigma_{A,B}$-orbit of $p_0$
such that $|p_i|$ is stricty increasing with $i$.
\qed

\begin{cor}\label{7.sols.cor}
  Let $S_{A,B}:=S_{a_1, a_2, b_1, b_2}$ be as in (\ref{m-s-3.thm.plan}.1).
  If $\max\{|a_2|+|a_1|, |b_2|+|b_1|\}\geq 7$, then there is a
  trivial solution  $p_0=(\pm 1, \pm 1, z_0)$ of a companion equation
  such that the  $\Sigma_{A,B}$-orbit of $p_0$ is infinite. 
\end{cor}

Proof.  Working with a suitable companion equation, we may assume  that
$$
|a_2|\leq |a_1|\qtq{and}  |a_2|+|a_1|\leq |b_2|+|b_1|=:d\geq 7.
$$
Choosing $y_0=\pm 1$ such that $b_2y_0^2$ and $b_1y_0$ have the same sign, we get that
$|B(y_0)|\geq |b_2|+|b_1|-2=d-2$. Thus 
$(x_1, y_1):=\sigma_x(x_0, y_0)$ such that $|x_1|\geq d-2$ and  $y_1=\pm 1$.
Then applying $\sigma_y$ we get $(x_2, y_2):=\sigma_y(x_1, y_1)$ where
$$
\begin{array}{rcl}
|y_2|&\geq & (d-2)^3-|a_2|(d-2)^2-|a_1|(d-2)-1\\
&\geq & (d-2)^3-\rdown{d/2}(d-2)^2-d(d-2)-1.
\end{array}
$$
It is easy to check that this is $>d+2$ for $d\geq 7$.
Thus the    $\Sigma_{A,B}$-orbit of $p_2$ is infinite
by Corollary~\ref{inf.orb.cor.1}. \qed

\begin{prop}\label{fin.triv.orb.prop}  
  For  $S_{A,B}:=S_{a_1, a_2, b_1, b_2}$  as in (\ref{m-s-3.thm.plan}.1),  the
  $\aut(S_{A,B})$-orbit of  every trivial solution  is finite iff $S_{A,B}$ is a  companion surface of one of the following.
  $$
  \begin{array}{lccl}
   (\ref{fin.triv.orb.prop}.1)\qquad \qquad  &    xyz-x^3-y^3-1 & = & -x^2-y^2,   \\
    (\ref{fin.triv.orb.prop}.2) &  xyz-x^3-y^3-1 & = & -2x^2-x-2y^2-y,\\
    (\ref{fin.triv.orb.prop}.3) &  xyz-x^3-y^3-1 & = & -2x^2-x-y^2,\\
   (\ref{fin.triv.orb.prop}.4) & xyz-x^3-y^3-1 & = & -x^2-2y^2-y.
  \end{array}
  $$
  (The last 2 are isomorphic to each other.)
\end{prop}

Proof.  By \cite[Thm.29]{k-vp-1}, $\langle \sigma_{A,B}\rangle$ has finite index in
$\aut_{\c}(S_{A,B})$.  Thus we need to enumerate all cases
for which the
$\Sigma_{A,B}$-orbit of  every trivial solution  is finite.
By Corollary~\ref{7.sols.cor}, this can happen only for
$|a_1|+|a_2|, |b_1|+ |b_2|\leq 6$. So we have a finite list
of possible exceptions.
The  computer program  is available at \url{https://web.math.princeton.edu/~jl5270/Mordell-Schinzel-code}. \qed

\begin{rem} There are 4 polynomials $t^3+a_2t^2+a_1t+1\in \z[t]$
  that map  $\{1, -1\}$ to itself.  These are
  $t^3-t+1,  t^3-t^2+1, t^3-t^2-2t+1, t^3-2t^2-t+1$.
  Thus, combining with  the list (\ref{fin.triv.orb.prop}.1--4), we see that either $\Sigma_{A,B}$ maps trivial solutions to trivial solutions, or the  $\Sigma_{A,B}$-orbit of trivial solutions is infinite.

  We do not have a direct proof of this.
  \end{rem}

\begin{say}[End of the proof of Theorem~\ref{m-s-3.thm}]
By Proposition~\ref{fin.triv.orb.prop},
there is a   trivial solution whose
$\aut(S_{A,B})$-orbit  is infinite, except for the 4 cases listed in
(\ref{fin.triv.orb.prop}.1--4).  

To handle these remaining cases, 
by Corollary~\ref{inf.orb.cor.1}, it is  enough to find for each one an
integral point $(x_0, y_0, z_0)$ with $\max\{|x_0|, |y_0|\}\geq 6$.
This was accomplished by a direct search. The smallest such solutions are
$$
  \begin{array}{lc}
   (\ref{fin.triv.orb.prop}.1)\qquad \qquad  &    (-7, -17, -47), \\
    (\ref{fin.triv.orb.prop}.2) & (293, -601, 1095),\\
    (\ref{fin.triv.orb.prop}.3) & (11, -13, 9),\\
   (\ref{fin.triv.orb.prop}.4) &  (-13, 11, 9).
  \end{array}
  $$
  The  computer program  is available at

  \url{https://web.math.princeton.edu/~jl5270/Mordell-Schinzel-code}. \qed
\end{say}

 \begin{say}[The geometry of the surfaces $S_{A,B}$]\label{geom.say}
  Let $\bar S_{A,B}\subset \p^3$ denote the closure of the affine cubic surface
  $S_{A,B}\subset \a^3$. The curve at infinity is the nodal cubic curve
  $(xyz=ax^3+by^3)\subset \p^2$, and   $\bar S_{A,B}$ has an $A_2$-singularity at the node. That is, in suitable local analytic coordinates, it can be written as $(uv=w^3)$.

  Conversely, assume that we are  over a field $k$ of characteristic different from 2 or 3. Let $S\subset \a^3$ be a cubic surface with projective closure
  $\bar S\subset \p^3$. Assume that the  curve at infinity is a nodal cubic, and  $\bar S$ has an $A_2$-singularity at the node. Then, in suitable coordinates,
  $S=S_{A, B}$ for some $A, B$ as in (\ref{not.1}.2);
   see \cite[Sec.2]{k-vp-1} for details.

  The presence of the $A_2$-singularity at infinity leads to the
  infinite automorphism group of such surfaces; see \cite{k-lk3,k-vp-1} for descriptions of such automorphism groups.
  If we allow $c=0$, the singularity becomes $A_3$-type,   $(uv=w^4)$.
  \end{say}

\begin{say}[Comments on \cite{MR51852} and \cite{MR3320468}]\label{ms.comments}
  From the point of view  of \cite{k-vp-1}, the proof in
  \cite{MR3320468} has 2 main steps, though they are not separated there.

A preliminary observation is  that if $|abc|>1$, then    the equation (\ref{not.1}.1)  has infinitely many $\z[(abc)^{-1}]$-integral solutions: we can take  $x, y$ to be arbitrary monomials in $a,b,c$, and set
$z=G(x,y)/(xy)$.
Call these {\it monomial solutions}   (though $z$ is not assumed monomial).
The automorphism   $\sigma_{A,B}$  of (\ref{g.isom.say}.2) is also defined over $\z[(abc)^{-1}]$.

The first step is then to show that for every $r>0$
 there is a
monomial point  $p_r:=(x_r, y_r, z_r)$ such that
$\sigma_{A,B}^r(p_r)$ is a $\z$-integral point.
The $x_r, y_r$ are explicitly given in the form
$a^{\lambda_r}b^{\mu_r}c^{\nu_r}$, where
$\lambda_r, \mu_r, \nu_r$ satisfy certain Fibonacci-type recursions;
see \cite[Formula~30]{MR3320468}.
This is surprising, but the proof  is relatively short; see  \cite[Lem.14]{MR3320468}.
This is clearly the plan in
\cite{MR51852}, though the general form in \cite{MR3320468} is a  nontrivial extension.

The second step is to show that these $\sigma_{A,B}^r(p_r)$ do give
infinitely many different points. This is clear in the examples
$xyz=ax^3+by^3+c$ considered in \cite{MR51852}, provided  $|abc|\geq 2$.
However,  neither \cite{MR51852} nor \cite[Sec.30]{MR249355} gives any hint on how the general case should be treated.

This is the harder part of  Schinzel's papers.
The treatments  require   subtle  inductive assumptions and estimates. At the end, for $m,n\geq 3$ this approach works whenever  $|abc|>1$; see \cite[Thm.4]{MR3320468}.

Thus we are left with the cases $a,b,c\in \{1, -1\}$.
The disadvantage is that inverting $abc$ does not give new solutions,
so we have only the
 trivial solutions   $(x_0, y_0, z_0)$ where $x_0, y_0\in \{1, -1\}$ to start with.
The advantage is that, by \cite[Thm.3]{k-vp-1} the automorphism group 
of $S_{A,B}$ is infinite, and in most cases the orbit of trivial solutions is infinite.

Our proof works for $m=n=3$  since there are only  finitely many
 exceptions, and they can be treated by hand.

\end{say}

\begin{say}[Counting solutions]
  Let $N_{A,B}(C)$ denote the number of integral points on $S_{A,B}$  satisfying
  $-C\leq x,y\leq C$. 
 Extensions of the  Batyrev--Manin conjectures to K3 surfaces
(see for example \cite{vanluijk}) 
 suggest the following.
 \begin{enumerate}
 \item Maps $\a^1\to S_{A,B}$ defined over $\z$ contribute 
    $O(C^{\epsilon})$  to $N_{A,B}(C)$.
 \item All other integral points contribute $(\log C)^r$  to $N_{A,B}(C)$ for some $r\in \n$.
 \end{enumerate}
 A detailed  study of integral points on cubic surfaces is given in 
 \cite{browning2024integralpointscubicsurfaces}.
 More limited computations of $N_{A,B}(C)$ were done in connection with \cite{k-vp-1}.  Both of these support the above general description, but unfortunately suggest different powers of $\log C$.

 It is not easy to find  maps $\a^1\to S_{A,B}$ defined over $\z$; see
 \cite{MR4458255}, \cite[Sec.5.5]{MR4823864}  or \cite[51]{k-vp-1} for some  examples and further references.
\end{say}

 \begin{ack} 
          We thank    B.~Grechuk, P.~Sarnak  and Sz.~Tengely  for    helpful comments and    references.

  Partial  financial support  to JK   was provided  by  the NSF under grant number
DMS-1901855 and by the Simons Foundation   under grant number SFI-MPS-MOV-00006719-02.
        \end{ack}


 \def\cprime{$'$} \def\cprime{$'$} \def\cprime{$'$} \def\cprime{$'$}
  \def\cprime{$'$} \def\dbar{\leavevmode\hbox to 0pt{\hskip.2ex
  \accent"16\hss}d} \def\cprime{$'$} \def\cprime{$'$}
  \def\polhk#1{\setbox0=\hbox{#1}{\ooalign{\hidewidth
  \lower1.5ex\hbox{`}\hidewidth\crcr\unhbox0}}} \def\cprime{$'$}
  \def\cprime{$'$} \def\cprime{$'$} \def\cprime{$'$}
  \def\polhk#1{\setbox0=\hbox{#1}{\ooalign{\hidewidth
  \lower1.5ex\hbox{`}\hidewidth\crcr\unhbox0}}} \def\cdprime{$''$}
  \def\cprime{$'$} \def\cprime{$'$} \def\cprime{$'$} \def\cprime{$'$}
\providecommand{\bysame}{\leavevmode\hbox to3em{\hrulefill}\thinspace}
\providecommand{\MR}{\relax\ifhmode\unskip\space\fi MR }
\providecommand{\MRhref}[2]{%
  \href{http://www.ams.org/mathscinet-getitem?mr=#1}{#2}
}
\providecommand{\href}[2]{#2}

 \bigskip

  Princeton University, Princeton NJ 08544-1000 

  \email{kollar@math.princeton.edu}
  
  \email{jenniferli@math.princeton.edu}
  
\end{document}